\documentclass[11pt,twoside,a4paper]{article}
\usepackage{amssymb}
\usepackage{amsfonts}
\usepackage{amsmath}
\usepackage{amsthm}
\usepackage{makeidx}
\usepackage[dvips]{graphicx}
\newtheorem{thm}{Theorem}

\newtheorem{lem}[thm]{Lemma}

\newcommand{\vertdot}{\circle*{0.4}}
\newcommand{\cirdot}{\circle*{0.3}}
\newcommand{\origo}{\circle{0,4}}
\newcommand{\Ve}{\mathcal{V}}
\newcommand{\Z}{\mathbb{Z}}

\newcommand{\R}{\mathbb{R}}
\newcommand{\indre}{\mathrm{int}}
\newcommand{\conv}{\mathrm{conv}}
\newcommand{\hyper}[2]{H(#1,#2)}
\newcommand{\nv}[2]{n(#1,#2)}
\newcommand{\nf}[2]{N(#1,#2)}
\newcommand{\vol}[1]{\mathrm{vol}(#1)}
\newcommand{\pair}[2]{\langle #1,#2 \rangle}

\title{Classification of terminal simplicial reflexive $d$-polytopes with $3d-1$ vertices.}
\author{Mikkel \O bro}
\def\address#1#2{\begingroup
\noindent\parbox[t]{7.8cm}{%
\small{\scshape\ignorespaces#1}\par\vskip1ex
\noindent\small{\itshape E-mail address}%
\/: #2\par\vskip4ex}\hfill%
\endgroup}%

\begin{document}

\maketitle

\begin{abstract}

We classify terminal simplicial reflexive $d$-polytopes with $3d-1$
vertices. They turn out to be smooth Fano $d$-polytopes. When $d$ is even there is 1 such polytope up to isomorphism, while there are 2 when $d$ is uneven.

\end{abstract}

\section{Introduction}

Let $N\cong \Z^d$ be a $d$-dimensional lattice, and let
$N_\R=N\otimes_{\Z}\R \cong \R^d$. Let $M$ be the dual lattice of $N$
and $M_\R$ the dual of $N_\R$. A \emph{reflexive} $d$-polytope $P$ in
$N_\R$ is a fully-dimensional convex lattice polytope, such that the
origin is contained in the interior and such that the dual polytope
$P^*:=\{x\in M_\R| \pair{x}{y}\leq 1 \ \forall y\in P\}$ is also a
lattice polytope. The notion of a reflexive polytope was introduced in
\cite{batyrev94}. Two reflexive polytopes
$P$ and $Q$ are called \emph{isomorphic} if there exists a bijective
linear map $\varphi:N_\R\to N_\R$, such that $\varphi(N)=N$ and
$\varphi(P)=Q$. For every $d\geq 1$ there are finitely many
isomorphism classes of reflexive $d$-polytopes, and for $d\leq 4$ they have been completely classified using computer algorithms (\cite{ks98},\cite{ks00}).

Simplicial reflexive $d$-polytopes have at most $3d$ vertices
(\cite{casagrande06} theorem 1). This upper bound is attained if and
only if $d$ is even and $P$ splits into $\frac{d}{2}$ copies of
\emph{del Pezzo 2-polytopes} $V_2=\conv\{\pm e_1,\pm e_2,\pm
(e_1-e_2)\}$, where $\{e_1,e_2\}$ is a basis of a $2$-dimensional lattice.

A reflexive polytope $P$ is called \emph{terminal}, if $N\cap P=0\cup \Ve(P)$. An important subclass of terminal simplicial reflexive polytopes is the class of \emph{smooth} reflexive polytopes, also known as \emph{smooth Fano} polytopes: A reflexive polytope $P$ is called \emph{smooth} if the vertices of every face $F$ of $P$ is a part of a basis of the lattice $N$. Smooth Fano $d$-polytopes have been intensively studied and completely classified up to dimension 4 (\cite{batyrev82},\cite{batyrev99},\cite{sato00},\cite{ww82}). In higher dimensions not much is known. There are classification results valid in any dimension, when the polytopes have few vertices (\cite{batyrev91},\cite{klein88}) or if one assumes some extra symmetries (\cite{casagrande03}, \cite{ewald88}, \cite{vk85}). Some of these results have been generalized to simplicial reflexive polytopes (\cite{nill05b}).

In this paper we classify terminal simplicial reflexive $d$-polytopes with $3d-1$ vertices for arbitrary $d$. It turns out that these are in fact smooth Fano $d$-polytopes.

\begin{thm} Let $P\subset N_\R$ be a terminal simplicial reflexive $d$-polytope with $3d-1$ vertices. Let $e_1,\ldots,e_d$ be a basis of the lattice $N$.

If $d$ is even, then $P$ is isomorphic to the convex hull of the points 
\begin{equation}
\label{polytope1}
\begin{array}{c}
e_1\ ,\ \pm e_2\ ,\ \ldots \ ,\ \pm e_d
\\
\pm(e_1-e_2)\ ,\ \ldots \ ,\ \pm(e_{d-1}-e_d).
\end{array}
\end{equation}
If $d$ is uneven, then $P$ is isomorphic to either the convex hull of the points
\begin{equation}
\label{polytope2}
\begin{array}{c}
\pm e_1\ ,\ \ldots \ ,\ \pm e_{d-1}\ ,\ e_d
\\
\pm(e_1-e_2)\ ,\ \ldots \ ,\ \pm(e_{d-2}-e_{d-1})\ ,\ e_1-e_d.
\end{array}
\end{equation}
or the convex hull of the points
\begin{equation}
\label{polytope3}
\begin{array}{c}
\pm e_1\ ,\ \ldots \ ,\ \pm e_d
\\
\pm (e_2-e_3)\ ,\ \ldots \ ,\ \pm (e_{d-1}-e_d).
\end{array}
\end{equation}
In particular, $P$ is a smooth Fano $d$-polytope.
\label{3d-1thm}
\end{thm}

A key concept in this paper is the notion of a special facet: A facet
$F$ of a simplicial reflexive $d$-polytope $P$ is called \emph{special}, if
the sum of the vertices $\Ve(P)$ of $P$ is a non-negative linear
combination of vertices of $F$. In particular, $\pair{u_F}{\sum_{v\in\Ve(P)}
v}\geq 0$, where $u_F\in M_\R$ is the unique element determined by $\pair{u_F}{F}=\{1\}$. The polytope $P$ is reflexive, so $\pair{u_F}{v}$ is an integer for
every $v\in\Ve(P)$. As $\pair{u_F}{v}\leq 1$ with equality if and only if
$v\in F$, there are at most $d$ vertices $v$ such that $\pair{u_F}{v}\leq
-1$. For simplicity, let $\hyper{F}{i}:=\{x\in N|\pair{u_F}{x}=i\}$, $i\in\Z$. It is well-known that at most $d$ vertices of $P$ are situated in $\hyper{F}{0}$ for any facet $F$ of $P$ (\cite{debarre} section 2.3 remarks 5(2)). If $P$ has $3d-1$ vertices and $F$ is a special facet of $P$, then
$$
d-1\leq|\Ve(P)\cap \hyper{F}{0}|\leq d,
$$
and there are only three possibilities for the placement of the $3d-1$ vertices of $P$ in the hyperplanes $\hyper{F}{i}$ as shown in the table below.
\begin{center}
\begin{tabular}{|c|ccc|}
\hline
 & Case 1 & Case 2 & Case 3\\
\hline
$|\Ve(P)\cap \hyper{F}{1}|$ & $d$ & $d$ & $d$\\
$|\Ve(P)\cap \hyper{F}{0}|$ & $d$ & $d$ & $d-1$\\
$|\Ve(P)\cap \hyper{F}{-1}|$ & $d-1$ & $d-2$ & $d$\\
$|\Ve(P)\cap \hyper{F}{-2}|$ & 0 & 1 & 0 \\
\hline
$|\Ve(P)|$ & $3d-1$ & $3d-1$ & $3d-1$\\
\hline
\end{tabular}
\end{center}
We prove theorem \ref{3d-1thm} by considering these three cases separately for terminal simplicial reflexive $d$-polytopes.

The paper is organised as follows: In section \ref{notation} we define some notation and prove some well-known basic facts about simplicial reflexive polytopes. In section \ref{afsnit3} we define the notion of special facets. In section \ref{afsnit4} we prove some lemmas needed in section \ref{bevis} for the proof of theorem \ref{3d-1thm}.
\\

\textbf{Acknowledgments}. The author would like to thank his advisor Johan P. Hansen for advice and encouragement.

\section{Notation and basic results}
\label{notation}

In this section we fix the notation and prove some basic facts about simplicial reflexive $d$-polytopes.

From now on $N$ denotes a $d$-dimensional lattice, $N\cong \Z^d$, $d\geq 1$, and $M$ denotes the dual lattice. Let $N_\R:=N\otimes_{\Z} \R$ and let $M_\R$ denote the dual of $N_\R$.

By $\conv K$ we denote the convex hull of a set $K$. A polytope is the convex hull of finitely many points, and a $k$-polytope is a polytope of dimension $k$. Recall that faces of a polytope of dimension 0 are called \emph{vertices}, while codimension 1 and 2 faces are called \emph{facets} and \emph{ridges}, respectively. The set of vertices of any polytope $P$ is denoted by $\Ve(P)$.

\subsection{Simplicial polytopes containing the origin in the interior}
\label{afsnit1}
A $d$-polytope $P$ in $N_\R$ is called \emph{simplicial} if every face of $P$ is a simplex.

\emph{In this section $P$ will be a simplicial $d$-polytope in $N_\R$ with $0\in\indre P$.}

For any facet $F$ of $P$, we define $u_F$ to be the unique
element in $M_\R$ where $\pair{u_F}{x}=1$ for every point $x\in F$. Certainly for any vertex $v$ and any facet $F$ of $P$, $\pair{u_F}{v}\leq 1$ with equality if and only if $v$ is a vertex of $F$.

We also define some points $u_F^v\in M_\R$ for any facet $F$ of $P$
and any vertex $v\in \Ve(F)$: $u_F^v$ is the unique element where $\pair{u_F^v}{v}=1$ and $\pair{u_F^v}{w}=0$ for every $w\in \Ve(F)\setminus\{v\}$. In other words, $\{u_F^v|v\in\Ve(F)\}$ is the basis of $M_\R$ dual to the
basis $\Ve(F)$ of $N_\R$.

When $F$ is a facet of $P$ and $v\in\Ve(F)$, there is a unique ridge
$R=\conv(\Ve(F)\setminus \{v\})$ of $P$ and a unique facet
$F'$ of $P$, such that $F\cap F'=R$. We denote this facet by
$\nf{F}{v}$ and call it a \emph{neighboring facet of $F$}. The set
$\Ve(\nf{F}{v})$ consists of the vertices $\Ve(R)$ of the ridge $R$
and a unique vertex $v'$, which we call a \emph{neighboring vertex of
  $F$} and denote it by $\nv{F}{v}$. See figure 1.

\begin{figure}
\begin{center}
\begin{picture}(12,7)
\thicklines
\put(1,2){\line(3,4){3}}
\put(1,2){\line(3,1){6}}
\put(4,6){\line(3,-2){3}}
\put(4,6){\line(6,-1){6}}
\put(7,4){\line(3,1){3}}
\put(1,2){\line(4,-1){4}}
\put(7,4){\line(-2,-3){2}}
\put(7,4){\cirdot}
\put(4,6){\cirdot}
\put(10,5){\cirdot}
\put(1,2){\cirdot}
\put(5,1){\cirdot}
\put(0.6,1.6){$v$}
\put(3.8,6.3){$w$}
\put(4,4){$F$}
\put(10.2,5){$\nv{F}{v}$}
\put(3,1.8){$\nf{F}{w}$}
\put(5,0.5){$\nv{F}{w}$}
\put(6.3,4.7){$\nf{F}{v}$}
\end{picture}
\caption{This illustrates the concepts of neighboring facets and
  neighboring vertices.}
\end{center}
\end{figure}
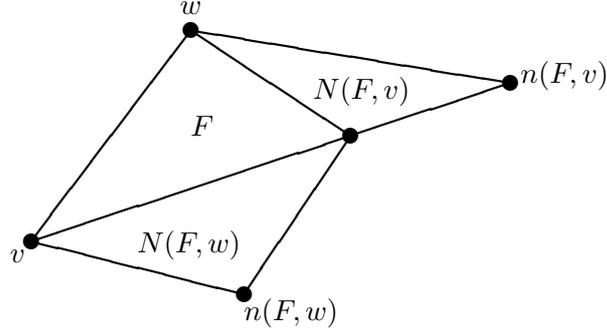

The next lemma shows how $u_F$ and $u_{F'}$ are related, when $F'$ is
a neighboring facet of the facet $F$.
\begin{lem}
\label{simple_uf_lemma}
Let $P\subset N_\R$ be a simplicial $d$-polytope containing the origin in the interior. Let $F$ be a facet of $P$ and $v\in \Ve(F)$. Let $F'$ be the neighboring facet $\nf{F}{v}$ and $v'$ the neighboring vertex $\nv{F}{v}$.

Then for any point $x\in N_\R$,
$$
\pair{u_{F'}}{x}=\pair{u_F}{x}+(\pair{u_{F'}}{v}-1)\pair{u_F^v}{x}.
$$
In particular,
\begin{itemize}
\item $\pair{u_F^v}{x}<0$ iff $\pair{u_{F'}}{x}>\pair{u_F}{x}$.
\item $\pair{u_F^v}{x}>0$ iff $\pair{u_{F'}}{x}<\pair{u_F}{x}$.
\item $\pair{u_F^v}{x}=0$ iff $\pair{u_{F'}}{x}=\pair{u_F}{x}$.
\end{itemize}
\begin{proof}
The vertices of $F$ span $N_\R$, and
$$
x=\sum_{w\in \Ve(F)} \pair{u_F^w}{x} w\ \ \ \textnormal{and}\ \ \ \pair{u_F}{x}=\sum_{w\in\Ve(F)} \pair{u_F^w}{x}.
$$
The vertices of the neighboring facet $F'=\nf{F}{v}$ are $\{v'\} \cup \Ve(F)\setminus \{v\}$. So
\begin{eqnarray*}
\pair{u_{F'}}{x}&=&\pair{u_F^v}{x} \pair{u_{F'}}{v} + \pair{u_F}{x}-\pair{u_F^v}{x}\\
&=&\pair{u_F}{x}+(\pair{u_{F'}}{v}-1)\pair{u_F^v}{x}.
\end{eqnarray*}
The vertex $v$ is not on the facet $F'$, and then the term $\pair{u_{F'}}{v}-1$ is negative. From this the equivalences follow.
\end{proof}
\end{lem}

\subsection{Simplicial reflexive polytopes}
\label{afsnit2}
A polytope $P\subset N_\R$ is called a \emph{lattice} polytope if
$\Ve(P)\subset N$. A lattice polytope is called \emph{reflexive}, if
$0\in\indre P$ and $\Ve(P^*)\subset M$, where
$$
P^*:=\{x\in M_\R\ |\ \pair{x}{y}\leq 1 \ \forall y\in P\}
$$
is the dual of $P$.

\emph{From now on $P$ denotes a simplicial reflexive $d$-polytope.}

Reflexivity guarentees that $u_F\in M$ for every facet $F$ of $P$, and every vertex of $P$ lies in one of the lattice hyperplanes
$$
\hyper{F}{i}:=\{x\in N\ |\ \pair{u_F}{x}=i\}\ \ \ ,\ i\in \{1,0,-1,-2,\ldots\}
$$
In particular, for every facet $F$ and every vertex $v$ of $P$: $v\notin F$ iff $\pair{u_F}{v}\leq 0$. This can put some restrictions on the points of $P$.

\begin{lem}
\label{coef_lemma}
Let $P$ be a simplicial reflexive polytope. For every facet $F$ of $P$ and every vertex $v\in\Ve(F)$ we have
$$
\pair{u_F}{x}-1\leq \pair{u_F^v}{x}
$$
for any $x\in P$. In case of equality, $x$ is on the facet $\nf{F}{v}$.
\begin{proof}
The inequality is obvious, when $\pair{u_F^v}{x}> 0$. So assume $\pair{u_F^v}{x}\leq
0$. Let $F'$ be the neighboring facet $\nf{F}{v}$.

Since $x\in P$, $\pair{u_{F'}}{x}\leq 1$ with equality iff $x\in F'$. From lemma \ref{simple_uf_lemma} we then have
$$
\pair{u_F}{x}-1\leq (1-\pair{u_{F'}}{v})\pair{u_F^v}{x}\leq \pair{u_F^v}{x}
$$
as $\pair{u_{F'}}{v}\leq 0$.
\end{proof}
\end{lem}

The next lemma concerns an important property of simplicial reflexive polytopes.

\begin{lem}[\cite{debarre} section 2.3 remarks 5(2), \cite{nill05} lemma 5.5]
Let $F$ be a facet and $x\in\hyper{F}{0}$ be
  vertex of a simplicial reflexive polytope $P$. Then $x$ is a neighboring vertex
  of $F$.

More precisely, for every $w\in\Ve(F)$ where $\pair{u_F^w}{x}<0$, $x$ is
  equal to $\nv{F}{w}$.

In particular, for every $w\in \Ve(F)$ there is at most one vertex
$x\in \hyper{F}{0}\cap \Ve(P)$, with $\pair{u_F^w}{x}<0$.

As a consequence, there are at most $d$ vertices of $P$ in $\hyper{F}{0}$.
\label{lemma1}

\begin{proof} Since $\pair{u_F}{x}=\sum_{w\in\Ve(F)}\pair{u_F^w}{x}=0$ and $x\neq 0$, there is at least
one $w\in\Ve(F)$ for which $\pair{u_F^w}{x}<0$. Choose such a $w$ and consider the neighboring facet $F'=\nf{F}{w}$.
By lemma \ref{simple_uf_lemma} we get that
$0<\pair{u_{F'}}{x}\leq 1$. As $P$ is reflexive, $\pair{u_{F'}}{x}=1$ and then $x=\nv{F}{w}$.

The remaining statements follow immediately.
\end{proof}
\end{lem}

\section{Special facets}
\label{afsnit3}
Now we define the notion of special facets, which will be of great use to us in the proof of theorem \ref{3d-1thm}.

\emph{$P$ is a simplicial reflexive $d$-polytope in this section.}

Consider the sum of all the vertices of $P$,
$$
\nu_P:=\sum_{v\in\Ve(P)} v.
$$
There exists at least one facet $F$ of $P$ such that $\nu_P$ is a non-negative linear combination of vertices of $F$, i.e. $\pair{u_F^w}{\nu_P}\geq 0$ for every $w\in\Ve(F)$. We call facets with this property \emph{special}.

Let $F$ be a special facet of $P$. In particular we have that
$$
0\leq \pair{u_F}{\nu_P},
$$
which implies that
\begin{equation}
\label{hsum}
0\leq \sum_{v\in\Ve(P)} \pair{u_F}{v}=\sum_{i\leq 1} i|\hyper{F}{i}\cap \Ve(P)|=d+\sum_{i\leq -1} i|\hyper{F}{i}\cap \Ve(P)|. 
\end{equation}
As there are at most $d$ vertices in $\hyper{F}{0}$ we can easily see
that $|\Ve(P)|\leq 3d$, which was first proved by Casagrande using a
similar argument (\cite{casagrande06} theorem 1). Notice that $\pair{u_F}{v}\geq -d$ for every vertex $v$ of $P$. Notice also, that when $|\Ve(P)|$ is close to $3d$, the vertices of $P$ tend to be in the hyperplanes $\hyper{F}{i}$ for $i\in\{1,0,-1\}$.

\section{Many vertices in $\hyper{F}{0}$}
\label{afsnit4}
We now study some cases of many vertices in $\hyper{F}{0}$, where $F$ is a facet of a simplicial reflexive $d$-polytope. The lemmas proven here are ingredients in the proof of theorem \ref{3d-1thm}.

\begin{lem}
\label{lemma2}
Let $F$ be a facet of a simplicial reflexive $d$-polytope $P$. Suppose there are at least $d-1$ vertices $v_1,\ldots,v_{d-1}$ in $\Ve(F)$, such that $\nv{F}{v_i}\in\hyper{F}{0}$ and $\pair{u_F^{v_i}}{\nv{F}{v_i}}=-1$ for every $1\leq i\leq d-1$.

Then $\Ve(F)$ is a basis of the lattice $N$.
\begin{proof}
Follows from statement 3 in \cite{nill05} lemma 5.5.
\end{proof}
\end{lem}

\begin{lem}
\label{basis_lem}
Let $P$ be a simplicial reflexive $d$-polytope, such that
$$
|\Ve(P)\cap \hyper{F}{0}|\geq d-1
$$
for every facet $F$ of $P$. Then there exists a facet $G$ of $P$, such that $\Ve(G)$ is a basis of $N$.
\begin{proof}
By lemma \ref{lemma2} we are done, if there exists a facet $G$ such that the set
$$
\{v\in\Ve(G)\ |\ \nv{G}{v}\in\hyper{G}{0}\ \textrm{and}\ \pair{u_G^v}{\nv{G}{v}}=-1\}
$$
is of size at least $d-1$. So we suppose that no such facet exists.

For every facet $F$ we denote the volume of the $d$-simplex
$\conv(\{0\}\cup \Ve(F))$ by $\vol{F}$. When $v_1,\ldots,v_d$ are the
vertices of $F$, the volume $\vol{F}$ is equal to $\frac{1}{d!}| \mathrm{det}A_F |$, where $A_F$ is the matrix
$$
A_F:= \left( \begin{array}{c}
v_1 \\
\vdots \\
v_d
\end{array}
\right) .
$$
The volume $\vol{\nf{F}{v_i}}$ of the neighboring facet $\nf{F}{v_i}$ is then
$$
\vol{\nf{F}{v_i}}=\frac{1}{d!}| \mathrm{det}A_{F'} |=\frac{|\pair{u_F^{v_i}}{\nv{F}{v_i}}|}{d!}|\mathrm{det} A_F|.
$$
Now, let $F_0$ be an arbitrary facet of $P$. There must be at least one vertex $v$ of $F_0$, such that $v'=\nv{F_0}{v}\in\hyper{F_0}{0}$, but $\pair{u_{F_0}^v}{v'}\neq -1$. Then $0>\pair{u_F^v}{v'}>-1$ by lemma \ref{coef_lemma}. Let $F_1$ denote the neighboring facet $\nf{F_0}{v}$. Then $\vol{F_0}>\vol{F_1}$.

We can proceed in this way to produce an infinite sequence of facets
$$
F_0,F_1,F_2,\ldots\ \ \ \textrm{where}\ \ \ \vol{F_0}>\vol{F_1}>\vol{F_2}>\ldots .
$$
But there are only finitely many facets of $P$. A contradiction.
\end{proof}
\end{lem}

\begin{lem}
\label{smalllemma}
Let $F$ be a facet of a simplicial reflexive polytope $P\in N_\R$. Let $v_1,v_2\in\Ve(F)$, $v_1\neq v_2$, and set $y_1=\nv{F}{v_1}$ and $y_2=\nv{F}{v_2}$. Suppose $y_1\neq y_2$, $y_1,y_2\in\hyper{F}{0}$ and $\pair{u_F^{v_1}}{y_1}=\pair{u_F^{v_2}}{y_2}=-1$.

Then there are no vertex $x\in \Ve(P)$ in $\hyper{F}{-1}$ with $\pair{u_F^{v_1}}{x}=\pair{u_F^{v_2}}{x}=-1$.
\begin{proof}
Suppose the statement is not true.

For simplicity, let $G=\conv(\Ve(F)\setminus \{v_1,v_2\})$. The vertex $x$ written as a linear combination of $\Ve(F)$ is then
$$
x=-v_1-v_2+\sum_{w\in\Ve(G)} \pair{u_F^w}{x} w.
$$
The vertices of the facet $F_1=\nf{F}{v_1}$ are $\{y_1\}\cup(\Ve(F)\setminus \{v_1\})$, where
$$
y_1=-v_1+\pair{u_F^{v_2}}{y_1}v_2+\sum_{w\in\Ve(G)} \pair{u_F^w}{y_1}w.
$$
In the basis (of $N_\R$) $F_1$ provides we have
$$
x=y_1+(-1-\pair{u_F^{v_2}}{y_1})v_2+\sum_{w\in\Ve(G)} \pair{u_F^w}{x-y_1} w
$$
The vertex $x$ is in $\hyper{F_1}{0}$ by lemma \ref{simple_uf_lemma}. Certainly, $\pair{u_F^{v_2}}{y_1}\leq 0$, otherwise we would have a contradiction to lemma \ref{coef_lemma}. On the other hand, $\pair{u_F^{v_2}}{y_1}\geq 0$, as $\nv{F}{v_2}\neq y_1$. So $\pair{u_F^{v_2}}{y_1}=0$ and $x=\nv{F_1}{v_2}$.

Similarly, $\pair{u_F^{v_1}}{y_2}=0$.
$$
y_2=-v_2+\sum_{w\in\Ve(G)} \pair{u_F^w}{y_2}w.
$$
But then $y_2$ and $x$ are both in $\hyper{F_1}{0}$ and both have negative $v_2$-coordinate. This is a contradiction to lemma \ref{lemma1}.
\end{proof}
\end{lem}

\subsection{The terminal case}

If we assume that the simplicial reflexive $d$-polytope $P$ is terminal, we can sharpen our results in case of $d$ vertices in $\hyper{F}{0}$ for some facet $F$ of $P$. Recall, that a reflexive polytope is called \emph{terminal} if $\Ve(P)\cup \{0\}=P\cap N$.

\begin{lem}
\label{lemh02}
Let $P$ be a terminal simplicial reflexive $d$-polytope. If there are $d$ vertices of $P$ in $\hyper{F}{0}$ for some facet $F$ of $P$, then 
$$
\Ve(P)\cap\hyper{F}{0}=\{ -y+z_y \ |\ y\in \Ve(F)\}
$$
where $z_y\in\Ve(F)$ for every $y$.

In particular $\Ve(F)$ is a basis of the lattice $N$.
\begin{proof} Let $y\in\Ve(F)$. By lemma \ref{lemma1} there exists exactly one vertex $x\in\hyper{F}{0}$, such that $x=\nv{F}{y}$. Conversely, there are no vertex $y'\neq y$ of $F$, such that $x=\nv{F}{y'}$. So $x$ is on the form
$$
x=-by+a_1w_1+\ldots+a_kw_k\ \ \ ,\ 0<b\leq 1\ ,\ 0<a_i\ \textnormal{and}\ w_i\in\Ve(F)\setminus\{y\}\ \forall i,
$$
where $b=\sum_{i=1}^k a_i$.

Suppose there exists a facet $G$ containing both $x$ and $y$. Then
$$
1+b=\pair{u_G}{x+by}=\pair{u_G}{a_1w_1+\ldots+a_kw_k}\leq \sum_{i=1}^k a_i=b.
$$
Which is a contradiction. So there are no such facets.

Consider the lattice point $z_y=x+y$. For any facet $G$ of $P$, $\pair{u_G}{z_y}\leq 1$ as both $\pair{u_G}{x},\pair{u_G}{y}\leq 1$ and both cannot be equal to 1. So $z_y$ is a lattice point in $P$. Since $P$ is terminal, $z_y$ is either a vertex of $P$ or the origin.

As $1=\pair{u_F}{x+y}=\pair{u_F}{z_y}$, $z_y$ must be a vertex of $F$ and $y\neq
z_y$. And then we're done.

The vertex set $\Ve(F)$ is a basis of $N$ by lemma \ref{lemma2}.
\end{proof}
\end{lem}

The proof of lemma \ref{lemh02} is inspired by proposition 3.1 in \cite{nill05}.

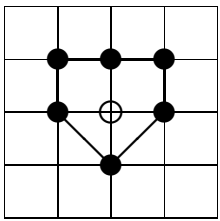
\begin{figure}
\begin{center}
\begin{picture}(4,4)
\linethickness{0.05mm}
\multiput(0,0)(1,0){5}{\line(0,1){4}}
\multiput(0,0)(0,1){5}{\line(1,0){4}}
\put(2,0){\line(0,1){4}}
\put(0,2){\line(1,0){4}}
\put(2,3){\vertdot}
\put(1,3){\vertdot}
\put(1,2){\vertdot}
\put(2,1){\vertdot}
\put(3,3){\vertdot}
\put(3,2){\vertdot}
\thicklines
\put(2,2){\origo}
\put(2,1){\line(-1,1){1}}
\put(1,2){\line(0,1){1}}
\put(1,3){\line(1,0){2}}
\put(3,3){\line(0,-1){1}}
\put(3,2){\line(-1,-1){1}}
\end{picture}
\caption{Terminality is important in lemma \ref{lemh02}: This is a simplicial reflexive (self-dual) $2$-polytope with 5 vertices. Consider the facet $F$ containing 3 lattice points. The two vertices in $\hyper{F}{0}$ is not on the form $-y+z_y$ for vertices $y,z_y\in \Ve(F)$.}
\end{center}
\end{figure}

\begin{lem} Let $F$ be a facet of a terminal simplicial reflexive $d$-polytope $P\subset N_\R$, such that $|\hyper{F}{0}\cap \Ve(P)|=d$. If $x\in\hyper{F}{-1}\cap P$, then $-x\in\Ve(F)$.
\label{smalllemma2}
\begin{proof}
The vertex set $\Ve(F)$ is a basis of the lattice $N$, and every vertex in $\hyper{F}{0}$ is of the form $-y+z$ for some $y,z\in\Ve(F)$ (lemma \ref{lemh02}).

Let $x$ be vertex of $P$ in $\hyper{F}{-1}$.
$$
x=\sum_{w\in\Ve(F)} \pair{u_F^w}{x} w,
$$
where $\pair{u_F^w}{x}\in\Z$ for every $w\in\Ve(F)$. If $\pair{u_F^w}{x}\leq -2$ for some $w\in\Ve(F)$, then $x=\nv{F}{w}$ (lemma \ref{coef_lemma}), which is not the case. So $\pair{u_F^w}{x}\geq -1$ for every $w\in\Ve(F)$. Furthermore, by lemma \ref{smalllemma} $x$ is only allowed one negative coordinate with respect to the basis $\Ve(F)$. The only possibility is then $x=-w$, where $w\in \Ve(F)$.
\end{proof}
\end{lem}

\section{Proof of main result}
\label{bevis}
In this section we will prove theorem \ref{3d-1thm}.

\emph{Throughout the section $P$ is a terminal simplicial reflexive
  $d$-polytope in $N_\R$ with $3d-1$ vertices, whose sum is $\nu_P$,
and $\{e_1,\ldots,e_d\}$ is a basis of the lattice $N$.}
$$
\nu_P:=\sum_{v\in\Ve(P)} v.
$$
By the existing classification we can check that theorem \ref{3d-1thm} holds for $d\leq 2$ (\cite{nill05} proposition 4.1). So we may assume that $d\geq 3$.

Let $F$ be a special facet of $P$, i.e. $\pair{u_F^w}{\nu_P}\geq 0$ for every $w\in\Ve(F)$. Of course, there are $d$ vertices of $P$ in $\hyper{F}{1}$. The remaining $2d-1$ vertices are in the hyperplanes $\hyper{F}{i}$ for $i\in\{0,-1,-2,\ldots,-d\}$, such that
$$
0\leq \pair{u_F}{\nu_P}=d+\sum_{i\leq -1} i\cdot |\Ve(P)\cap \hyper{F}{i}|.
$$
So there are three cases to consider.
\begin{center}
\begin{tabular}{|c|ccc|}
\hline
 & Case 1 & Case 2 & Case 3\\
\hline
$|\Ve(P)\cap \hyper{F}{1}|$ & $d$ & $d$ & $d$\\
$|\Ve(P)\cap \hyper{F}{0}|$ & $d$ & $d$ & $d-1$\\
$|\Ve(P)\cap \hyper{F}{-1}|$ & $d-1$ & $d-2$ & $d$\\
$|\Ve(P)\cap \hyper{F}{-2}|$ & 0 & 1 & 0 \\
\hline
$|\Ve(P)|$ & $3d-1$ & $3d-1$ & $3d-1$\\
\hline
\end{tabular}
\end{center}
We will consider these cases seperately.

\begin{description}
\item[Case 1.] There are $d$ vertices in $\hyper{F}{0}$, so by lemma \ref{lemh02} $\Ve(F)$ is a basis of $N$. We may then assume that $\Ve(F)=\{e_1,\ldots,e_d\}$.

The sum of the vertices is a lattice point on $F$, since $\pair{u_F}{\nu_P}=1$. As $P$ is terminal, this must be a vertex $e_i$ of $F$, say $\nu_P=e_1$. Then a facet $F'$ of $P$ is a special facet iff $e_1\in \Ve(F')$.

There are $d-1$ vertices in $\hyper{F}{-1}$, so by lemma \ref{smalllemma2} we get
$$
\Ve(P)\cap \hyper{F}{-1}=\{-e_1,\ldots,-e_{j-1},-e_{j+1},\ldots,-e_d\},
$$
for some $1\leq j\leq d$. Now, there are two possibilities: $j=1$ or $j\neq 1$, that is $-e_1\notin\Ve(P)$ or $e_1\in\Ve(P)$.

\begin{description}
\item[$-e_1\notin\Ve(P)$.] Then $-e_i\in\Ve(P)$ for every $2\leq
i\leq d$. There are $d$ vertices in $\hyper{F}{0}$, so by lemma
\ref{lemh02} there is a vertex $-e_1+e_{a_1}$, which we can assume to be
$-e_1+e_2$.

Consider the facet $F'=\nf{F}{e_2}$. This is a special
facet, so we can show that
$$
\Ve(P)\cap \hyper{F'}{-1}=\Ve(-F')\setminus \{-e_1\}.
$$
The vertex $-e_1+e_2$ is in the hyperplane $\hyper{F'}{-1}$. So $e_1-e_2$ is a vertex of $F'$ (lemma \ref{smalllemma2}), and then of $P$.

For every $3\leq i\leq d$ we use
the same procedure to show that $-e_i+e_{a_i}$ and $-e_{a_i}+e_i$ are
vertices of $P$. This shows that $d$ is even and that $P$ is isomorphic to the convex hull of the points in (\ref{polytope1}).

\item[$-e_1\in\Ve(P)$.] We may assume $-e_d\notin \Ve(P)$. The sum of the vertices $\Ve(P)$ is $e_1$, so there are exactly two vertices in $\hyper{F}{0}$ of the form $-e_k+e_1$ and $-e_l+e_1$, $k\neq l$. We wish to show that $k=d$ or $l=d$. This is obvious for $d=3$. So suppose $d\geq 4$ and $k,l\neq d$, that is $-e_k,-e_l\in\Ve(P)$.

Consider the facet $F'=\nf{F}{e_k}$, which is a special facet. So by the arguments above we get that
$$
\Ve(P)\cap \hyper{F'}{-1}=\Ve(-F')\setminus \{-e_d\}.
$$
As $\Ve(F')=\{e_1,\ldots,e_{k-1},e_{k+1},e_d,-e_k+e_1\}$, we have that $-e_1+e_k$ must be a vertex of $P$.

In a similar way we get that $-e_1+e_l$ is a vertex of $P$. But this is a contradiction. So $k$ or $l$ is equal to $d$, and without loss of generality, we can assume that $k=2$ and $l=d$.

For $3\leq i\leq d-1$ we proceed in a similar way to get that both $-e_i+e_{a_i}$ and $-e_{a_i}+e_i$ are vertices of $P$.

And so we have showed that $d$ must be uneven and that $P$ is isomorphic to the convex hull of the points in (\ref{polytope2}).
\end{description}

\item[Case 2.] Since $\pair{u_F}{\nu_P}=0$, the sum of the vertices is the origin, so
  every facet of $P$ is special. There are $d$ vertices in $\hyper{F}{0}$, so $\Ve(F)$ is a basis of $N$ (lemma \ref{lemh02}). Without loss of generality, we can assume $\Ve(F)=\{e_1,\ldots,e_d\}$. By lemma \ref{smalllemma2}
$$
x\in \Ve(P)\cap \hyper{F}{-1}\ \ \Longrightarrow \ \ x=-e_i\ \textnormal{for some $1\leq i\leq d$.}
$$
Consider the single vertex $v$
  in the hyperplane $\hyper{F}{-2}$. If $\pair{u_F^{e_j}}{v}>0$ for some $j$ then
  $\pair{u_{F'}}{v}<-2$ for the facet $F'=\nf{F}{e_j}$ (lemma
  \ref{simple_uf_lemma}), which is not the case as $F'$ is special. So
  $\pair{u_F^{e_j}}{v}\leq 0$ for every $1\leq j\leq d$. As $v$ is a primitive lattice
  point we can without loss of generality assume $v=-e_1-e_2$.

There are $d$ vertices in $\hyper{F}{0}$, so there is a vertex of the form $-e_1+e_j$ for some
$j\neq 1$. If $j=2$, then $-e_1\in\conv\{-e_1+e_2,-e_1-e_2\}$ which is not the case as $P$ is terminal. So we may assume $j=3$. In $\hyper{F}{0}$ we also find the vertex $-e_2+e_i$ for some
$i$. A similar argument yields $i\neq 1$.

Let $G=\nf{F}{e_1}$. Then $\Ve(G)$ is a basis of the lattice $N$. Write $v$ in this basis.
$$
v=(-e_1+e_3)-e_3-e_2.
$$
As $i\neq 1$, $-e_2+e_i$ is in $\hyper{G}{0}$ and is equal to $\nv{G}{e_2}$ (lemma \ref{lemma1}).

If $v\neq \nv{G}{e_3}$ we have a contradiction to lemma \ref{smalllemma}. Therefore $v=\nv{G}{e_3}$, and $\conv\{v,-e_1+e_3,e_2\}$ is a face of $P$.

As $e_3$ and $-e_1+e_3$ are vertices of $P$, there are at least two
vertices of $P$ with positive $e_3$-coordinate (with respect to the basis $F$
provides). There is exactly one vertex in $\hyper{F}{0}$ with negative $e_3$-coordinate, namely $-e_3+e_k$ for some $k$. Any other has to be in $\hyper{F}{-1}$. The vertices of $P$ add to 0, so the point $-e_3$ must be a vertex of $P$.

But $-e_3=-(-e_1+e_3)+v+e_2$, which cannot be the case as $P$ is simplicial.

We conclude that case 2 is not possible.

\item[Case 3.] In this case we also have $\pair{u_F}{\nu_P}=0$, so every
  facet is special. Case 2 was not possible, so $-1\leq u_G(v)\leq 1$
  for any facet $G$ and any vertex $v$ of $P$. By lemma
  \ref{basis_lem} we may assume that $\Ve(F)$ is a basis of $N$, say $\Ve(F)=\{e_1,\ldots,e_d\}$. As case 2 was not possible, 
$$
\Ve(P)\cap \hyper{F}{-1}=\{-e_1,\ldots,-e_d\}.
$$
Let $x$ be any vertex in $\hyper{F}{0}$. There exists at least one $1\leq i\leq d$, such that $\pair{u_F^{e_i}}{x}=-1$. Consider the facet $F'=\nf{F}{e_i}$: $\Ve(F')=(\Ve(F)\setminus \{e_i\}) \cup \{x\}$ is a basis of $N$. Case 2 was not possible, so $\Ve(P)\cap \hyper{F'}{-1}=-\Ve(F')$, which implies $-x\in\Ve(P)$.

This shows that $P$ is centrally symmetric and $d$ must be uneven. By \cite{nill05} theorem 5.9 $P$ is isomorphic to the convex hull of the points in (\ref{polytope3}).

\end{description}

This ends the proof of theorem \ref{3d-1thm}.

\address{Department of Mathematical Sciences \\
University of \AA rhus \\
8000 \AA rhus C \\
Denmark
}
{oebro@imf.au.dk}

\end{document}